\documentclass[review]{elsarticle}
\usepackage[nottoc]{tocbibind}
\usepackage[ruled,vlined,onelanguage]{algorithm2e}
\newtheorem{lemma}{Lemma}
\newtheorem{remark}{Remark}
\newtheorem{theorem}{Theorem}
\newtheorem{proposition}{Proposition}
\usepackage{lineno,hyperref}
\modulolinenumbers[5]
\usepackage[utf8]{inputenc}
\usepackage{amsmath}
\usepackage{graphicx}
\usepackage{float}
\usepackage{subcaption}
\modulolinenumbers[1]
\usepackage{siunitx}
\usepackage{amsmath}
\usepackage[capitalise]{cleveref}
\usepackage{siunitx}
\journal{arXiv}

\usepackage{amsfonts}

\begin{document}

\begin{frontmatter}
\title{The transcritical Bogdanov Takens bifurcation with boundary due to the risk perception on a recruitment epidemiological model} 



\author[1]{Jocelyn A. Castro-Echeverr\'ia}
\author[2]{Fernando Verduzco}
\author[3]{Jorge X. Velasco-Hern\'andez}

\address[1]{Tecnologico de Monterrey, México / Instituto de Matem\'aticas UNAM Unidad Juriquilla,México}
\address[2]{Universidad de Sonora, México}
\address[3]{Instituto de Matem\'aticas UNAM Unidad Juriquilla,México}

\begin{abstract}
We analyze an epidemiological model with treatment and recruitment considering the risk perception. In this model, we consider an exponential function as a recruitment rate. We have found that this model undergoes the transcritical Bogdanov-Takens bifurcation with boundary, where the system experiences the transcritical bifurcation between the disease-free equilibrium point and the endemic equilibrium point. The Hopf bifurcation also arises at the endemic equilibrium point, this is, the appearance or disappearance of a limit cycle, and finally, the Homoclinic bifurcation which transforms the limit cycle into a homoclinic cycle, starting and ending at the disease-free equilibrium point.
\end{abstract}

\end{frontmatter}


\section{Introduction}
Behavior is a factor that pervades social interactions in all animal species. Particularly, transmission processes are strongly affected by behavioral decisions of individuals \cite{Mena2020,Banks2017,Ngonghala2015,Just2016}. During the COVID-19 pandemic, the collective community behavior had a huge impact on disease spread. This factor was one of the main reasons why many mathematical models faild in forecasting particular epidemic events such as peak dates, epidemic size, reproductive number and so forth \cite{Bavel2020,Sepulveda2021,Thoren2024}. Population behaviour was affected, not only by the control and mitigation measures instated by the all governments, but also by the risk perception of the community. Risk averse and risk prone behaviors vary within individuals and affect their collective behavior \cite{Nishimi2022,Chan2020,Espinoza2021}. Nowadays, frontier research in epidemiological dynamics involves the development of methodologies to measure the impact of behaviour on outbreaks and eradication of disease \cite{Bedson2021}. How to incorporate behavioral factors into useful mathematical models is also a main challenge in the postcovid era. 
There are several studies that have investigated the role of behavioral components in disease transmission. For example, in \cite{levin_insights_2021}, it is analyzed the relationship between mobility across different cities, registered by smartphone GPS data, and the evolution of the COVID-19 disease. The effect of mitigation of the social risk perception due to control policies  has been studied in \cite{ye_trust_2020,usherwood_model_2021,jamieson_race_2021}. Sexually-transmitted diseases are prime examples of the role of behavior in epidemic outcomes. Sexual intercourse involves partner selection, a paradigmatic behavioral trait. Therefore, the relevance of behavior in the recruitment of sexual partners were notoriously evident during AIDS pandemic, directly affecting the duration of that pandemic \cite{Velasco1996}. Other studies on social evitation and immigration are documented in \cite{Brauer2008, brauer_spatial_2019, Meng2021, Zou2022}.



Behavior change has been reidentified as a main and important driver of epidemic events and a key factor for the control and mitigation of outbreaks \cite{levin_insights_2021,ye_trust_2020,usherwood_model_2021, jamieson_race_2021,kellerer_behavior_2021,ha_changes_2022}. Epidemiological models have incorporated behavioral aspects in the study of several diseases \cite{Brauer1993, Velasco1996, Brauer2005, Ajbar2021}.  As mentioned above, sexually-transmitted diseases, notably gonorrhea and AIDS epidemic, have undelined the importance of behavioral change in the spread of the disease \cite{hethcote_gonorrhea_1984,Velasco1996}. Clearly, sexual behavior involves contagion risks that  depend on gender, education, and socio-economic status \cite{fentahun_risky_2014, lucas_schooling_2019}. Behavioral components, in particular core group of behavior, was recognized as a fundamental aspect of gonorrhea transmission in the pioneering work of \cite{hethcote_gonorrhea_1984} and was later  extended to AIDS \cite{Brauer1993,Velasco1996}. To date, the global impact of the COVID-19 pandemic has underlined, once again, that demographic, economic, and behavioral factors are necessary for the more efficient design of public health policies \cite{bhattacharyya_modelling_2021,Bai2021,Ajbar2021,Meng2021,levin_insights_2021,Zou2022,ghosh_mathematical_2022}.

Behavioral change can be approximated through the effective contact rate \cite{Brauer2011} and, by the use of the core group idea \cite{hethcote_gonorrhea_1984} recruitment into the core group may depend on the perceived prevalence level within the population where the disease is spreading under a mass-action mixing law and, thus, a linear incidence rate \cite{Velasco1996}.  In this paper, we revisit an already published model for prevalence-dependent recruitment \cite{Velasco1996} and concentrate on the dynamics of the bifurcation reported in that work. This dynamical system experiences a very interesting phenomenon, the so-called Takens-Bogdanov (T-B) bifurcation that appears as the result of the interaction between the disease-free equilibrium point and a boundary equilibrium point, where the homoclinic curve arises. The T-B bifurcation has been shown previously in a long list of epidemic models, for instance in \cite{Wang2004,Alexander2005,Moghadas2006,Xiao2015,Li2015,Lu2020,Zhang2022} all of them considering a nonlinear incidence rate. In addition, in \cite{Shan2014, Misra2022}, it is considered the number of hospital beds as a compartment, generating also the T-B bifurcation.  This is evidence that this kind of dynamic behavior is present in real-life models not only for epidemics but also in some others biological models, as in \cite{Liu2016} where the T-B appears for a prey-depredator model. Our main interest is the analysis of the T-B bifurcation and the relationship with the risk perception factor.

\section{Model setup}
The model that we analyze in this work revisits the paper \cite{Velasco1996} and is given in Eq \ref{sistema-original}. The model assumes that there is a general population $P$ where a subpopulation or core group $N$, suffers a directly-transmitted infectious diseases. The model is a $SI$ Kermack-McKendrick model with treatment $U$. Originally, this model was presented as approximating the core group behavior observed during the AIDS pandemic in certain cities \cite{Velasco1996} and introduces the concept of prevalence-dependent recruitment where susceptible individuals $S$ recruit into the core group depending upon the level of prevalence observed in the subpopulation. This process is modelled using an exponential function of the form $\exp^{-\left(a_1 \frac{I}{N}+a_2 \frac{U}{N}\right)}$ where $a_1>0$ and $a_2<0$ are called the risk coefficients. This function increases or decreases dependent upon prevalence $I/N$ or the level of treatment $U/N$ in the core group. The model equations are:

\begin{eqnarray}\label{sistema-original}
\dot{P}&=&b\left(P+S\right)+\hat{b}\left(I+U\right)-\left[\theta\left(S,I,U\right)+\mu\right]P \nonumber\\
\dot{S}&=&\theta\left(S,I,U\right)P-\beta SI/N-\mu S+\gamma U \nonumber\\
\dot{I}&=&\left(\beta S+\hat{\beta} U\right)I/N-\left(\mu+\tau\right)I \nonumber\\
\dot{U}&=&\tau I-\hat{\beta}UI/N-\left(\mu+\gamma\right)U,
\end{eqnarray}
where $\hat{\beta}=\left(1-\eta\right)\beta$, $\theta\left(S,I,U\right)=\exp^{-\left(a_1 \frac{I}{N}+a_2 \frac{U}{N}\right)}$, with $$N=U+S+I \text{ y } T=P+N$$ 
for $P\geq0$, $S\geq0$, $I\geq0$ and $U\geq0$. The state variables are defined by

\begin{figure}[H]]\centering
                    \includegraphics[height=1.8in]{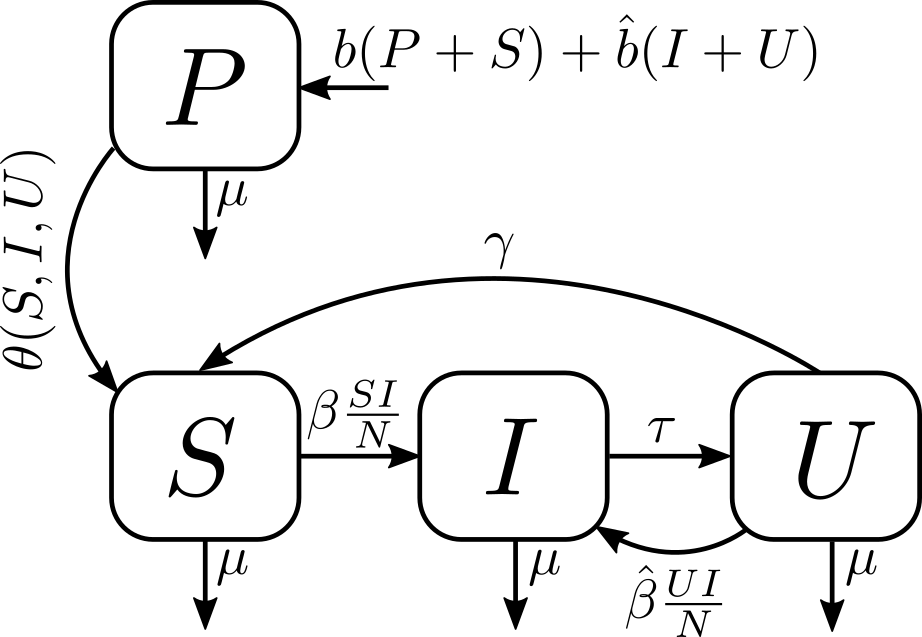}
                    \captionsetup{labelformat=empty}
                    \caption{Diagram for system (\ref{sistema-original}).}
                    \label{Fig-diagram1}
\end{figure}

\begin{table}[H]
\centering
\begin{tabular}{clc}
\hline
\textbf{State variables}     & \textbf{Description} & \textbf{Rank}           \\ \hline
$P$        &  Non-core group           & $\mathbb{R}^+$\\
$S$        &  Susceptible on core-group       & $\mathbb{R}^+$\\
$I$        &  Infectious on core-group           & $\mathbb{R}^+$\\
$U$        &  Treated population on core-group            & $\mathbb{R}^+$\\ 
$N$        &  Total population on core-group          & $\mathbb{R}^+$\\
$T$        &  Total population          & $\mathbb{R}^+$\\
\hline
\end{tabular}
\caption{State variables of model (\ref{sistema-original}).}\label{Tab:estados}
\end{table}
with parameters 
\begin{table}[H]
\centering
\begin{tabular}{cll}
\hline
\textbf{Parameters} & \textbf{Description}                    & \textbf{Rango} \\ \hline
$b$        &  Birth rate of susceptible population  &  [0,\,1] \\
$\hat{b}$  &  Birth rate of infected population  &  [0,\,1] \\
$\eta$     &  Treatment efficacy &  [0,\,1] \\         
$a_1$      &  Risk perception factor due to prevalence  & $\mathbb{R}^+\cup\{0\}$\\
$a_2$      &  Risk perception factor due to treatment & $\mathbb{R}^-\cup\{0\}$\\
$\mu$      &  Mortality rate                &  [0,\,1] \\
$\gamma$   &  Recovery rate                 &  [0,\,1]  \\
$\tau$     &  Tratment rate                &  [0,\,1] \\
$\beta$    &  Infectious rate               &  [0,\,1] \\ \hline
\end{tabular}
\caption{Parameters of model (\ref{sistema-original}).}\label{Tab:parametros}
\end{table}
System (\ref{sistema-original}) diagram is shown in Figure \ref{Fig-diagram1}.

Let us note that
\begin{eqnarray}
\dot{T}&=&\dot{P}+\dot{S}+\dot{I}+\dot{U}\nonumber\\
       &=&b\left(P+S\right)+\hat{b}\left(I+U\right)-\mu\left(P+S+I+U\right),\nonumber
\end{eqnarray}
Hence, if we consider the birth rates of the susceptible population $b$ and infected population $\hat{b}$ the same and equal to the mortality rate $\mu$, i.e. $b=\hat{b}=\mu$, we have that $\dot{T}=0$, this means that the total population $T$ remains constant over the time.\\

Now, under the assumption that $T$ is constant and given that $P\left(S,I,U\right)=T-S-I-U$ we are able to rewrite
\begin{eqnarray}
\dot{P}&=&\dot{T}-\dot{S}-\dot{I}-\dot{U}\nonumber\\
       &=&-\left(\dot{S}+\dot{I}+\dot{U}\right),\nonumber
\end{eqnarray}
this allows us to reduce the system (\ref{sistema-original}), from four to three dimensions, by eliminating the equation for $\dot{P}$. From this previous explanation, we obtain the following reduce system in $\mathbb{R}^3$
\begin{eqnarray}\label{sistema-reducido}
\dot{S}&=&\theta\left(S,I,U\right)\left(T-S-I-U\right)-\beta SI/N-\mu S+\gamma U \nonumber\\
\dot{I}&=&\left(\beta S+\hat{\beta} U\right)I/N-\left(\mu+\tau\right)I\nonumber\\
\dot{U}&=&(\tau - \hat{\beta} U/N) I - (\mu + \gamma) U
\end{eqnarray}
with $N=S+I+U$.

\subsection{Desease free equilibrium point}\label{section_dfe}
The disease-free equilibrium point $E_0=(S_0,I_0,U_0)$ is such that if $I_0=0$ then $\dot{I}=0$, and also from $\dot{U}=0$ we have $U_0=0$. Later, $\theta\left(S_0,0,0\right)=1$, consequently 
\begin{eqnarray}
\dot{S}=0 &\Leftrightarrow&T-S_0-\mu S_0=0 \nonumber\\
        &\Leftrightarrow&S_0=\frac{T}{\mu+1}. \nonumber\\
\end{eqnarray}

Therefore, the disease-free equilibrium point has coordinates
\begin{eqnarray}\label{dfe}
E_0=(S_0,I_0,U_0)=\left(\frac{T}{\mu+1},\;0,\;0\right).
\end{eqnarray}

\begin{remark} It is worth noting that the disease-free equilibrium point is always present, regardless of the specific parameter values.
\end{remark}

For $y=\left(I,S,U\right)^T$, according to the algorithm described in \cite{van2002reproduction}, to compute the basic reproductive number $R_0$ of system (\ref{sistema-reducido}) we might rewrite

\begin{eqnarray*}
\dot{y}=\mathcal{F}(y)-\mathcal{V}(y),
\end{eqnarray*}
where 
\begin{eqnarray*}
\mathcal{F}(y)&=&\left(\left(\beta S+\hat{\beta}U\right)I/N ,\; 0, \;0\right)^T,\\
\mathcal{V}^+(y)&=&\left(0,\; \theta(S,I,U)\; P(S,I,U)+\gamma U,\; \tau I\right)^T,\\
\mathcal{V}^-(y)&=&\left(\left(\mu+\tau\right)I,\; \beta SI/N+\mu S, \; \hat{\beta} UI/N+\left(\mu+\gamma\right)U\right)^T,
\end{eqnarray*}
with $\mathcal{V}(y)=\mathcal{V}^-(y)-\mathcal{V}^+(y)$.

As we only have one infected compartment we define
$$F=\frac{\partial F_1}{\partial I}\left(E_0\right)=\beta \text{ y } V=\frac{\partial V_1}{\partial I}\left(E_0\right)=\mu+\tau,$$
where $F_1(y)=\left(\beta S+\hat{\beta}U\right)I/N$ and $V_1(y)=\left(\mu+\tau\right)I$, which implies that the next generation matrix is the scalar number $K=FV^{-1}$, therefore,
\begin{eqnarray}\label{R0_def}
R_0=\dfrac{\beta}{\mu+\tau}.
\end{eqnarray}

Let $\dot{x}=f(x)$, with $x=\left(S,I,U\right)^T$, then the Jacobian matrix of system (\ref{sistema-reducido}) is given by
\begin{eqnarray}\label{jacobiana}
Df(x)=\left(\begin{array}{ccc}
\frac{\partial \theta}{\partial S} P-\theta-\frac{\beta I}{N}-\mu & \frac{\partial \theta}{\partial I} P-\theta-\frac{\beta S}{N} & \frac{\partial \theta}{\partial I} P-\theta+\gamma\\
\frac{\beta I}{N} &\frac{\beta S+\hat{\beta}U}{N}-\left(\mu+\tau\right) & \frac{\hat{\beta}I}{N}\\
0 &\tau-\frac{\hat{\beta}U}{N} & \frac{\hat{\beta}I}{N}-\left(\mu+\gamma\right)
\end{array}\right),
\end{eqnarray}
being $P=P(S,I,U)$.\\

Now, by evaluating matrix (\ref{jacobiana}) at the disease-free equilibrium point $E_0=(S_0,I_0,U_0)=\left(\frac{T}{\mu+1},0,0\right)$ we obtain
\begin{eqnarray}
Df(E_0)=\left(\begin{array}{ccc}
-\left(\mu+1\right) & -\left(\beta+a_1\mu+1\right) & \gamma-\left(a_2\mu+1\right)\\
0 & \beta-\left(\mu+\tau\right) & 0 \\
0 & \tau & -\left(\gamma+\mu\right)
\end{array}\right)\nonumber
\end{eqnarray}
whose eigenvalues are 
\begin{eqnarray*}
\lambda_1&=&-\left(\mu+1\right),\\
\lambda_2 &=&\beta-\left(\mu+\tau\right),\\
\lambda_3&=&-\left(\mu+\gamma\right).
\end{eqnarray*}

Note that both $\lambda_1$ and $\lambda_3$ are negative values. We can rewrite $\lambda_2$ as $$\lambda_2 =\left(R_0-1\right)\left(\mu+\tau\right).$$ Therefore, the stability of the disease-free equilibrium point depends on the basic reproductive number, $R_0$.

\begin{theorem}\label{Teo_dfe}
The point $E_0=(S_0,I_0,U_0)=\left(\frac{T}{\mu+1},\;0,\;0\right)$ is the desease-free equilibrium point of system (\ref{sistema-reducido}), this equilibrium is always present and its stability is given as follows:
\begin{itemize}
    \item If $R_0<1$ the disease-free equilibrium point is stable,
    \item If $R_0>1$ the disease-free equilibrium point is unstable,
\end{itemize}
where $R_0=\frac{\beta}{\mu+\tau}$ is the basic reproductive number.
\end{theorem}

\begin{figure}[H]\centering
                    \includegraphics[height=3in]{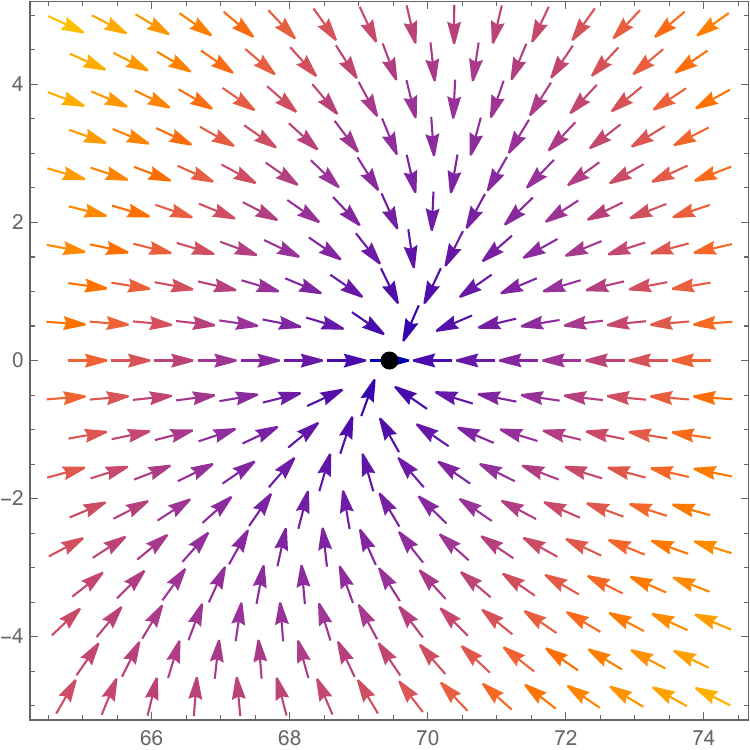}
                    \captionsetup{labelformat=empty}
                    \caption{Simulation of system (\ref{sistema-reducido}) at the invariant plane $I=0$.}
                    \label{Fig-planoinv}
\end{figure}

Whenever $I$ equals to zero, it follows that $\dot{I}$ is also zero. Furthermore, due to the signs of $\lambda_1$ and $\lambda_3$, we can deduce the following proposition:

\begin{proposition}
The plane $S$-$U$ remains invariant, and the point $(S_0,U_0)$ is an equilibrium point which is a stable node. Please refer to Figure \ref{Fig-planoinv}.
\end{proposition}

\subsection{Endemic equilibrium point}
Upon examining system (\ref{sistema-reducido}), we search for any equilibrium point $E_1=\left(S_1,I_1,U_1\right)$ such that $I_1\neq 0$. In this case, $N_1=S_1+I_1+U_1\neq0$, then, by taking $\dot{I}=0$, we have
\begin{eqnarray*}
\left(\beta S_1+\hat{\beta} U_1\right)-\left(\mu+\tau\right)N_1&=&0,\\
\left(\beta S_1+\hat{\beta} U_1\right)-\left(\mu+\tau\right)\left(S_1+I_1+U_1\right)&=&0,\\
S_1\left(\beta -\mu-\tau\right)+U_1\left(\hat{\beta}-\mu-\tau\right)-I_1 \left(\mu+\tau\right)&=&0,
\end{eqnarray*}
from where 
\begin{eqnarray}\label{eq_l1}
I_1\left(S_1,U_1\right)&=&S_1\left(\frac{\beta}{\mu+\tau}-1\right)+U_1\left(\frac{\beta\left(1-\eta\right)}{\mu+\tau}-1\right)\nonumber\\
&=&S_1\left(R_0-1\right)+U_1\left(R_0\left(1-\eta\right)-1\right)
\end{eqnarray}
Additionally, for the equilibrium point to exist in the biological sense, it is necessary that $I_1\left(S_1,U_1\right)>0$, this is
\begin{eqnarray*}
&& S_1\left(R_0-1\right)+U_1\left(R_0\left(1-\eta\right)-1\right)>0\\
&\Leftrightarrow& R_0\left(S_1+U_1\left(1-\eta\right)\right)-\left(S_1+U_1\right)>0, \\
&\Leftrightarrow& R_0>\frac{S_1+U_1}{S_1+U_1\left(1-\eta\right)},\\
&\Leftrightarrow& R_0>1,
\end{eqnarray*}
since $0<\eta<1$. In conclusion, the existence of the endemic equilibrium point requires $R_0>1$, indicating that it arises when the disease-free equilibrium is unstable.

Now, for the point $E_1=\left(S_1,I_1\left(S_1,U_1\right),U_1\right)$ according to the system (\ref{sistema-reducido}), $\dot{U}=0$ if and only if
\begin{eqnarray*}
&&(\tau - \hat{\beta} U_1) I_1\left(S_1,U_1\right) - (\mu + \gamma) U_1 N_1=0,\\
&&(\tau - \hat{\beta} U_1) I_1\left(S_1,U_1\right) - (\mu + \gamma) U_1 \left(S_1+I_1\left(S_1,U_1\right)+U_1\right)=0,
\end{eqnarray*}
which is an equation only in the $S_1$ and $U_1$ variables, whose solution is
\begin{eqnarray}\label{eq_s1}
    S_1(U_1)=C_0\;U_1 ,
\end{eqnarray}
where
\begin{eqnarray*}
    C_0=\frac{ \sqrt{D_0}+\gamma +\eta  \mu +\tau -R_0  (\eta -1) (\mu -\tau )}{2 \tau  (R_0-1)},
\end{eqnarray*}
with $$D_0=2 \tau  (\gamma  (2 \eta -1)+\eta  \mu +(\eta -1) R_0 (\mu +\tau ))+(\gamma +\eta  \mu -(\eta -1) R_0(\mu +\tau ))^2+\tau ^2.$$
    
Finally, by solving $\dot{S}=0$ for $E_1=\left(S_1(U_1),I_1(S_1(U_1),U_1),U_1\right)$ in system (\ref{sistema-reducido}), we obtain the expression for $U_1$ in terms of the parameters only:
\begin{eqnarray}\label{eq_u1}
  U_1= \frac{2 T(R_0-1) \tau  }{R_0 \left(\gamma +\eta  \mu +2 \eta  \tau +\sqrt{D_0}-(\eta -1) R_0 (\mu +\tau )-\tau \right) \left(\mu  e^{C1}+1\right)},
\end{eqnarray}
by taking
\begin{eqnarray*}
C_1=\frac{(R_0-1) \left(a_1 \left(\gamma +\eta  \mu +\sqrt{D_0}-\eta  R_0 (\mu +\tau )+R_0 (\mu +\tau )-\tau \right)+2 a_2 \tau \right)}{R_0 \left(\gamma +\eta  \mu +2 \eta  \tau +\sqrt{D_0}-\eta  R_0 (\mu +\tau )+R_0 (\mu +\tau )-\tau \right)}. 
\end{eqnarray*}

\begin{theorem}\label{Teo_eep}
If $R_0>1$, then, the point $E_1=\left(S_1(U_1),I_1(S_1(U_1),U_1),U_1\right)$ defined by equations (\ref{eq_l1}), (\ref{eq_s1}) and (\ref{eq_u1}) is an endemic equilibrium point of system (\ref{sistema-reducido}).
\end{theorem}





\section{Bifurcation analysis}

\subsection{Transcritical bifurcation (Forward bifurcation).}

As stated in theorem (\ref{Teo_dfe}), system (\ref{sistema-reducido}) always has a disease-free equilibrium point. However, the stability of this point varies based on the value of $R_0$. Additionally, when $R_0>1$, an endemic equilibrium point emerges while the disease-free equilibrium point becomes unstable.\\

Recall from equation(\ref{eq_l1}) that for the endemic equilibrium point $E_1$, we have
\begin{eqnarray*}
I_1\left(S_1,U_1\right)&=&S_1\left(R_0-1\right)+U_1\left(R_0\left(1-\eta\right)-1\right),
\end{eqnarray*}
If we assume $R_0>1$ and $R_0 \rightarrow 1$, then,
\begin{eqnarray*}
I_1\left(S_1,U_1\right)\rightarrow-\eta U_1,
\end{eqnarray*}
It is not possible for $I_1(S_1,U_1)$ to be negative, which means that $U_1$ must be equal to zero in order to have a feasible option. This leads to $I_1$ approaching zero, which aligns with the disease-free equilibrium point $E_0$, resulting in a transcritical bifurcation at the disease-free equilibrium point. As $E_0$ is located at the boundary of the system (\ref{sistema-reducido}), this bifurcation is referred to as a transcritical bifurcation with boundary, or more specifically, a forward bifurcation in this case.

\begin{theorem}\label{Teo_trans}
      The system (\ref{sistema-reducido}) undergoes a transcritical bifurcation with boundary (forward bifurcation) at the disease-free equilibrium point when $R_0\rightarrow 1$, as follows:
    \begin{itemize}
        \item[a)] If $R_0<1$ the disease-free equilibrium point $E_0$ is the only equilibrium of the system, and it is stable.
        \item[b)] If $R_0>1$ the disease-free equilibrium point becomes unstable, and an endemic equilibrium point emerges.
    \end{itemize}    
\end{theorem}

\subsection{Transcritical Bogdanov-Takens bifurcation}
\begin{figure}[H]\centering
                    \includegraphics[height=4in]{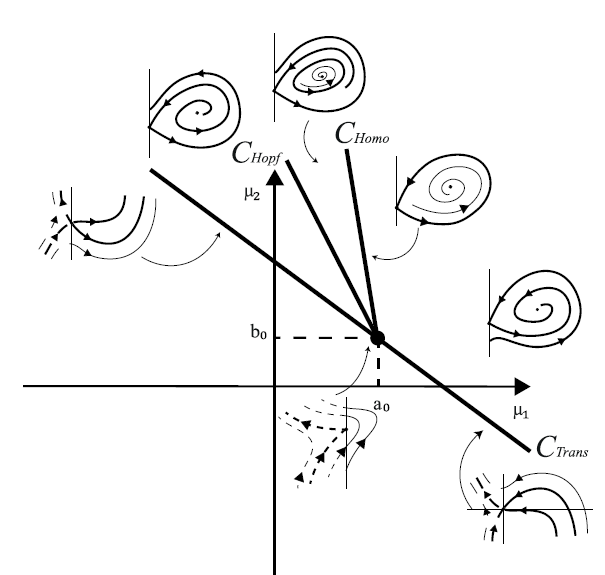}
                    \caption{tBt Bifurcation with boundary diagram presented in \cite{castro2020bifurcation}.}
                    \label{tbt-diagram}
\end{figure}

The transcritical Bogdanov-Takens bifurcation (tBT) is a variation of the typical Bogdanov-Takens bifurcation. The main difference is that, instead of the saddle-node bifurcation, the system undergoes the transcritical bifurcation, along with the Hopf and Homoclinic bifurcation. In the system (\ref{sistema-reducido}), there is evidence that this bifurcation occurs at the disease-free equilibrium point $E_0=(S_0,I_0,U_0)$, where both $I_0$ and $U_0$ are equal to zero. This equilibrium lies on the boundary of the system's domain since $S,I,U\geq0$. Therefore, we are dealing with a tBT bifurcation with a boundary. This phenomenon is well-described in \cite{castro2020bifurcation}, where the bifurcation diagram presented in Figure (\ref{tbt-diagram}) is shown.

\subsubsection{Center manifold analysis}
Consider the system given in (\ref{sistema-reducido}) re-writen as
\begin{eqnarray}\label{sistema_1}
    \dot{x}&=&f(x,\delta),
\end{eqnarray}
where $x=\left(S,I,U\right)^T$, $\delta=\left(\delta_1,\delta_2\right)^T$, and $\delta_1$, $\delta_2$ defined as
\begin{eqnarray}\label{deltas}
    \delta_1&=&R_0-1\nonumber\\
    \delta_2&=&\gamma+\mu.
\end{eqnarray}
The following Lemma is a direct result of the analysis in Section \ref{section_dfe}.
\begin{lemma}
For system (\ref{sistema_1}), for $E_0=\left(S_0,I_0,U_0\right)$ defined in (\ref{dfe}) and $\delta_0=\left(0,0\right)^T$, the following properties hold:
    \begin{itemize}
        \item[$(i)$] $f(E_0,\delta_0)=0$.
        \item[$(ii)$] $\sigma\left\{A=Df(E_0,\delta_0)\right\}=\left\{\lambda_2=0,\;\lambda_3=0 \text{ and } \lambda_0=\lambda_1<0 \right\}$.
    \end{itemize}
\end{lemma}
In this case, the matrix A is similar to the Jordan Canonical form of a nilpotent matrix in blocks:
\begin{eqnarray}\label{J0_matrix}
J=P^{-1}AP\left(\begin{array}{ccc}
0 & 1 & 0\\
0 & 0 & 0\\
0 & 0 & \lambda_0
\end{array}\right),
\end{eqnarray}
where matrix $P$ contains the eigenvectors of matrix $A$.\\

For a detailed analysis of system (\ref{sistema_1}), a change of variables is needed. First, we compute its Taylor Series expansion around the disease-free equilibrium point $E_0$:
\begin{eqnarray}
    \dot{x}&=&f(E_0,\delta_0)+Df(E_0,\delta_0)\;\left(x-E_0\right)+\frac{1}{2}\;D^2f(E_0,\delta_0)\;\left(x-E_0,x-E_0\right)\nonumber\\
    &&+\frac{\partial^2f}{\partial x \partial \delta }(E_0,\delta_0)\;\left(\delta,x-E_0\right)+\frac{1}{6}\;D^3f(E_0,\delta_0)\;\left(x-E_0,x-E_0,x-E_0\right)\nonumber\\
    &&+\frac{\partial^3f}{\partial^2 x \partial \delta }(E_0,\delta_0)\;\left(\delta,x-E_0,x-E_0\right)+\dots
\end{eqnarray}
Let $y=P^{-1}\left(x-E_0\right)$, then $Py=x-E_0$, therefore
\begin{eqnarray}\label{cambio1}
    \dot{y}&=&P^{-1}\;\dot{x},\nonumber\\
    &=&P^{-1}\left(APy+\frac{1}{2}\;D^2f(E_0,\delta_0)\;\left(Py,Py\right)+\frac{\partial^2f}{\partial x \partial \delta }(E_0,\delta_0)\;\left(\delta,Py\right)\right.\nonumber\\
    &&\left.+\frac{1}{6}\;D^3f(E_0,\delta_0)\;\left(Py,Py,Py\right)+\frac{\partial^3f}{\partial^2 x \partial \delta }(E_0,\delta_0)\;\left(\delta,Py,Py\right)+\dots\right),\nonumber\\
    &=&P^{-1}APy+\frac{1}{2}\;P^{-1}D^2f(E_0,\delta_0)\;\left(Py,Py\right)+P^{-1}\frac{\partial^2f}{\partial x \partial \delta }(E_0,\delta_0)\;\left(\delta,Py\right)\nonumber\\
    &&+\frac{1}{6}P^{-1}\;D^3f(E_0,\delta_0)\;\left(Py,Py,Py\right)+P^{-1}\frac{\partial^3f}{\partial^2 x \partial \delta }(E_0,\delta_0)\;\left(\delta,Py,Py\right)+\dots\nonumber\\
     &=&Jy+\frac{1}{2}\;P^{-1}D^2f(E_0,\delta_0)\;\left(Py,Py\right)+P^{-1}\frac{\partial^2f}{\partial x \partial \delta }(E_0,\delta_0)\;\left(\delta,Py\right)\nonumber\\
    &&+\frac{1}{6}P^{-1}\;D^3f(E_0,\delta_0)\;\left(Py,Py,Py\right)+P^{-1}\frac{\partial^3f}{\partial^2 x \partial \delta }(E_0,\delta_0)\;\left(\delta,Py,Py\right)+\dots\nonumber,\\
\end{eqnarray}
Now, considering $y=(w,z)^T$, $w\in\mathbb{R}^2$ and $z\in\mathbb{R}$, we can propose a center manifold
\begin{eqnarray}
    z=H(w,\delta)=\frac{1}{2}w^TH_1w+\delta^TH_2w+H_3(w)+\delta_1G_{21}(w)+\delta_2G_{22}(w),
\end{eqnarray}
it is necessary to this function to preserve the homological equation
\begin{eqnarray}
  \dot{z}&=&DH(w,\delta)\dot{w}\nonumber\\
    &=&\frac{\partial H}{\partial w_1}\dot{w}_1+\frac{\partial H}{\partial w_2}\dot{w}_2,
\end{eqnarray}
by solving this equation we can get the dynamical over the centre manifold.\\

After some computations made at Mathematica Software and by following the algorithm presented in \cite{kuznetsov2005practical} Appendix B, we conclude that the dynamical behavior over the center manifold, when $\delta=\delta_0$,has the form
\begin{eqnarray}
  \dot{w}_1&=&w_2\nonumber\\
  \dot{w}_1  &=&b_2w_1w_2+b_3w_1^2w_2+b_4w_1^3w_2+\ldots,
\end{eqnarray}
Moreover, it is possible to make $b_3$ and $b_4$ equal to zero by choosing appropriate eigenvectors of matrix $A$.

Unfortunately, according to \cite{kuznetsov2005practical}, this scenario has not yet been examined in current literature. It is also considered the most challenging case of this type of unfolding, and as a result, the analysis is beyond the scope of the existing tools available.

Although we are unable to provide proof of the existence of the tBT bifurcation with boundary, we do have numerical evidence indicating the emergence of a limit cycle of system (\ref{sistema-reducido}) around the endemic equilibrium point $E_1$, which corresponds to the Hopf Bifurcation. By setting $R_0$ close to 1, we can observe that the limit cycle is almost identical to a homoclinic cycle.

\section{Simulations}
In this section, we show some simulations as numerical evidence of the tBT bifurcation of system (\ref{sistema-reducido}).

\begin{figure}[H]\centering
                    \includegraphics[height=1.2in]{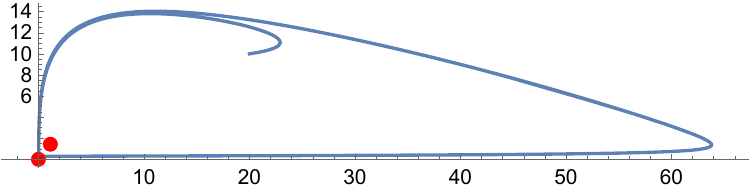}
                    \caption{Simulation at the $S$-$I$ plane, a transcritical bifurcation ocurrs between the disease-free and the endemic equilibrium points, meanwhile the cycle arises, converging to the homoclinic cycle.}
                    \label{Fig-simuhomo}
\end{figure}

Figure \ref{Fig-simuhomo} shows the phase portrait simulation for system (\ref{sistema-reducido}) generated under the following parameters values:
$$a_1=8, \;a_2=-2.3, \;\beta =1, \;\eta =0.41, \;\gamma =0.46, \;\mu =0.44, \;\tau =0.002, \;T=100.$$

\begin{figure}[H]\centering
                    \includegraphics[height=1.5in]{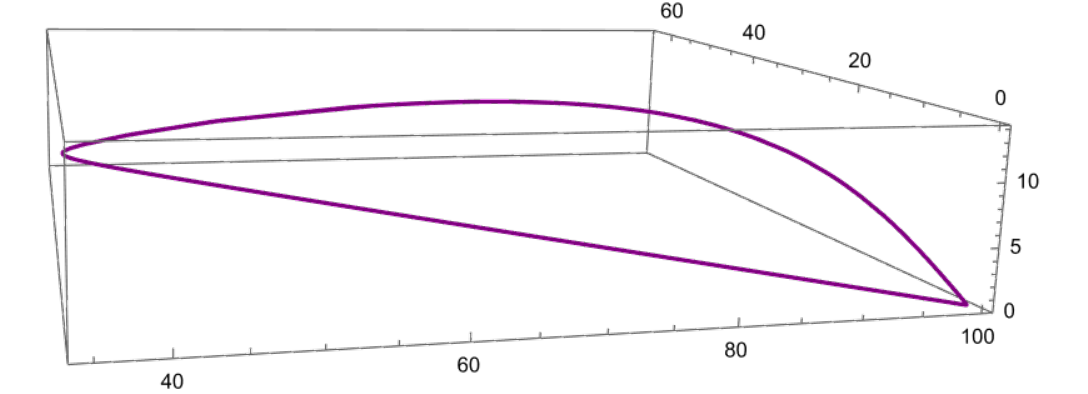}
                    \caption{3D simulation of the cycle around the endemic equlibrium converging to the homoclinic curve at the disease-free equilibrium.}
                    \label{Fig-simuhomo2}
\end{figure}

 
 \section*{Acknowledgemets}
\noindent JACE and JXVH acknowledge support from grants UNAM PAPIIT IV100220. 



\bibliography{mybibfile.bib}

\begin{thebibliography}{10}
\expandafter\ifx\csname url\endcsname\relax
  \def\url#1{\texttt{#1}}\fi
\expandafter\ifx\csname urlprefix\endcsname\relax\def\urlprefix{URL }\fi
\expandafter\ifx\csname href\endcsname\relax
  \def\href#1#2{#2} \def\path#1{#1}\fi

\bibitem{Mena2020}
R.~H. Mena, J.~X. Velasco-Hernandez, N.~B. Mantilla-Beniers, G.~A. Carranco-Sapiéns, L.~Benet, D.~Boyer, I.~P. Castillo, Using posterior predictive distributions to analyse epidemic models: Covid-19 in mexico city, Physical Biology 17 (11 2020).
\newblock \href {https://doi.org/10.1088/1478-3975/abb115} {\path{doi:10.1088/1478-3975/abb115}}.

\bibitem{Banks2017}
H.~T. Banks, K.~Bekele-Maxwell, R.~A. Everett, L.~Stephenson, S.~Shao, J.~Morgenstern, \href{http://link.springer.com/10.1007/s11538-017-0282-5}{Dynamic modeling of problem drinkers undergoing behavioral treatment}, Bulletin of Mathematical Biology 79 (2017) 1254--1273.
\newblock \href {https://doi.org/10.1007/s11538-017-0282-5} {\path{doi:10.1007/s11538-017-0282-5}}.
\newline\urlprefix\url{http://link.springer.com/10.1007/s11538-017-0282-5}

\bibitem{Ngonghala2015}
C.~N. Ngonghala, G.~D. Leo, M.~Pascual, D.~C. Keenan, A.~Dobson, M.~H. Bonds, General ecological models for human subsistence, health and poverty, Nature: ecology and evolution (2015) 1--21\href {https://doi.org/10.1038/s41559-017-0221-8} {\path{doi:10.1038/s41559-017-0221-8}}.

\bibitem{Just2016}
W.~Just, J.~Saldana, \href{http://arxiv.org/abs/1606.08788}{Oscillations in epidemic models with spread of awareness}, Journal of Mathematical Biology (2016).
\newblock \href {https://doi.org/10.1007/s00285-017-1166-x} {\path{doi:10.1007/s00285-017-1166-x}}.
\newline\urlprefix\url{http://arxiv.org/abs/1606.08788}

\bibitem{Bavel2020}
J.~J. Bavel, K.~Baicker, P.~S. Boggio, V.~Capraro, A.~Cichocka, M.~Cikara, M.~J. Crockett, A.~J. Crum, K.~M. Douglas, J.~N. Druckman, J.~Drury, O.~Dube, N.~Ellemers, E.~J. Finkel, J.~H. Fowler, M.~Gelfand, S.~Han, S.~A. Haslam, J.~Jetten, S.~Kitayama, D.~Mobbs, L.~E. Napper, D.~J. Packer, G.~Pennycook, E.~Peters, R.~E. Petty, D.~G. Rand, S.~D. Reicher, S.~Schnall, A.~Shariff, L.~J. Skitka, S.~S. Smith, C.~R. Sunstein, N.~Tabri, J.~A. Tucker, S.~van~der Linden, P.~van Lange, K.~A. Weeden, M.~J. Wohl, J.~Zaki, S.~R. Zion, R.~Willer, Using social and behavioural science to support covid-19 pandemic response (5 2020).
\newblock \href {https://doi.org/10.1038/s41562-020-0884-z} {\path{doi:10.1038/s41562-020-0884-z}}.

\bibitem{Sepulveda2021}
J.~Sepulveda, La respuesta de méxico al covid-19: Estudio de caso (2021).

\bibitem{Thoren2024}
H.~Thorén, P.~Gerlee, \href{https://royalsocietypublishing.org/doi/10.1098/rsos.230803}{Model uncertainty, the covid-19 pandemic, and the science-policy interface}, Royal Society Open Science 11 (2 2024).
\newblock \href {https://doi.org/10.1098/rsos.230803} {\path{doi:10.1098/rsos.230803}}.
\newline\urlprefix\url{https://royalsocietypublishing.org/doi/10.1098/rsos.230803}

\bibitem{Nishimi2022}
K.~Nishimi, B.~Borsari, B.~P. Marx, R.~C. Rosen, B.~E. Cohen, E.~Woodward, D.~Maven, P.~Tripp, A.~Jiha, J.~D. Woolley, T.~C. Neylan, A.~O'Donovan, \href{https://www.sciencedirect.com/science/article/pii/S2211335521003624}{Clusters of covid-19 protective and risky behaviors and their associations with pandemic, socio-demographic, and mental health factors in the united states}, Preventive Medicine Reports 25 (2022) 101671.
\newblock \href {https://doi.org/https://doi.org/10.1016/j.pmedr.2021.101671} {\path{doi:https://doi.org/10.1016/j.pmedr.2021.101671}}.
\newline\urlprefix\url{https://www.sciencedirect.com/science/article/pii/S2211335521003624}

\bibitem{Chan2020}
H.~F. Chan, A.~Skali, D.~A. Savage, D.~Stadelmann, B.~Torgler, Risk attitudes and human mobility during the covid-19 pandemic, Scientific Reports 10 (12 2020).
\newblock \href {https://doi.org/10.1038/s41598-020-76763-2} {\path{doi:10.1038/s41598-020-76763-2}}.

\bibitem{Espinoza2021}
B.~Espinoza, M.~Marathe, S.~Swarup, M.~Thakur, Asymptomatic individuals can increase the final epidemic size under adaptive human behavior, Scientific Reports 11 (12 2021).
\newblock \href {https://doi.org/10.1038/s41598-021-98999-2} {\path{doi:10.1038/s41598-021-98999-2}}.

\bibitem{Bedson2021}
J.~Bedson, L.~A. Skrip, D.~Pedi, S.~Abramowitz, S.~Carter, M.~F. Jalloh, S.~Funk, N.~Gobat, T.~Giles-Vernick, G.~Chowell, J.~R. de~Almeida, R.~Elessawi, S.~V. Scarpino, R.~A. Hammond, S.~Briand, J.~M. Epstein, L.~Hébert-Dufresne, B.~M. Althouse, A review and agenda for integrated disease models including social and behavioural factors, Nature Human Behaviour 5 (2021) 834--846.
\newblock \href {https://doi.org/10.1038/s41562-021-01136-2} {\path{doi:10.1038/s41562-021-01136-2}}.

\bibitem{levin_insights_2021}
R.~Levin, D.~L. Chao, E.~A. Wenger, J.~L. Proctor, \href{https://www.nature.com/articles/s43588-021-00125-9}{Insights into population behavior during the {COVID}-19 pandemic from cell phone mobility data and manifold learning}, Nature Computational Science 1~(9) (2021) 588--597.
\newblock \href {https://doi.org/10.1038/s43588-021-00125-9} {\path{doi:10.1038/s43588-021-00125-9}}.
\newline\urlprefix\url{https://www.nature.com/articles/s43588-021-00125-9}

\bibitem{ye_trust_2020}
M.~Ye, Z.~Lyu, \href{https://linkinghub.elsevier.com/retrieve/pii/S027795362030736X}{Trust, risk perception, and {COVID}-19 infections: {Evidence} from multilevel analyses of combined original dataset in {China}}, Social Science \& Medicine 265 (2020) 113517.
\newblock \href {https://doi.org/10.1016/j.socscimed.2020.113517} {\path{doi:10.1016/j.socscimed.2020.113517}}.
\newline\urlprefix\url{https://linkinghub.elsevier.com/retrieve/pii/S027795362030736X}

\bibitem{usherwood_model_2021}
T.~Usherwood, Z.~LaJoie, V.~Srivastava, \href{http://www.nature.com/articles/s41598-021-91514-7}{A model and predictions for {COVID}-19 considering population behavior and vaccination}, Scientific Reports 11~(1) (2021) 12051.
\newblock \href {https://doi.org/10.1038/s41598-021-91514-7} {\path{doi:10.1038/s41598-021-91514-7}}.
\newline\urlprefix\url{http://www.nature.com/articles/s41598-021-91514-7}

\bibitem{jamieson_race_2021}
T.~Jamieson, D.~Caldwell, B.~Gomez-Aguinaga, C.~Doña-Reveco, \href{https://www.mdpi.com/1660-4601/18/21/11113}{Race, {Ethnicity}, {Nativity} and {Perceptions} of {Health} {Risk} during the {COVID}-19 {Pandemic} in the {US}}, International Journal of Environmental Research and Public Health 18~(21) (2021) 11113.
\newblock \href {https://doi.org/10.3390/ijerph182111113} {\path{doi:10.3390/ijerph182111113}}.
\newline\urlprefix\url{https://www.mdpi.com/1660-4601/18/21/11113}

\bibitem{Velasco1996}
J.~X. Velasco-Hern{\'a}ndez, F.~Brauer, C.~Castillo-Chavez, Effects of treatment and prevalence-dependent recruitment on the dynamics of a fatal disease, Mathematical Medicine and Biology: A Journal of the IMA 13~(3) (1996) 175--192.

\bibitem{Brauer2008}
F.~Brauer, P.~van~den Driessche, L.~Wang, \href{https://doi.org/10.1016/j.mbs.2008.05.001}{Oscillations in a patchy environment disease model}, Mathematical Biosciences 215~(1) (2008) 1--10.
\newblock \href {https://doi.org/10.1016/j.mbs.2008.05.001} {\path{doi:10.1016/j.mbs.2008.05.001}}.
\newline\urlprefix\url{https://doi.org/10.1016/j.mbs.2008.05.001}

\bibitem{brauer_spatial_2019}
F.~Brauer, C.~Castillo-Chavez, Z.~Feng, \href{https://doi.org/10.1007/978-1-4939-9828-9_14}{Spatial {Structure} in {Disease} {Transmission} {Models}}, in: F.~Brauer, C.~Castillo-Chavez, Z.~Feng (Eds.), Mathematical {Models} in {Epidemiology}, Texts in {Applied} {Mathematics}, Springer, New York, NY, 2019, pp. 457--476.
\newblock \href {https://doi.org/10.1007/978-1-4939-9828-9_14} {\path{doi:10.1007/978-1-4939-9828-9_14}}.
\newline\urlprefix\url{https://doi.org/10.1007/978-1-4939-9828-9_14}

\bibitem{Meng2021}
L.~Meng, W.~Zhu, \href{https://doi.org/10.1155/2021/5401253}{Generalized {SEIR} epidemic model for {COVID}-19 in a multipatch environment}, Discrete Dynamics in Nature and Society 2021 (2021) 1--12.
\newblock \href {https://doi.org/10.1155/2021/5401253} {\path{doi:10.1155/2021/5401253}}.
\newline\urlprefix\url{https://doi.org/10.1155/2021/5401253}

\bibitem{Zou2022}
Y.~Zou, W.~Yang, J.~Lai, J.~Hou, W.~Lin, \href{https://doi.org/10.1007/s11538-021-00958-5}{Vaccination and quarantine effect on {COVID}-19 transmission dynamics incorporating chinese-spring-festival travel rush: Modeling and simulations}, Bulletin of Mathematical Biology 84~(2) (Jan. 2022).
\newblock \href {https://doi.org/10.1007/s11538-021-00958-5} {\path{doi:10.1007/s11538-021-00958-5}}.
\newline\urlprefix\url{https://doi.org/10.1007/s11538-021-00958-5}

\bibitem{kellerer_behavior_2021}
J.~D. Kellerer, M.~Rohringer, D.~Deufert, \href{https://onlinelibrary.wiley.com/doi/10.1111/phn.12918}{Behavior in the use of face masks in the context of {COVID}‐19}, Public Health Nursing 38~(5) (2021) 862--868.
\newblock \href {https://doi.org/10.1111/phn.12918} {\path{doi:10.1111/phn.12918}}.
\newline\urlprefix\url{https://onlinelibrary.wiley.com/doi/10.1111/phn.12918}

\bibitem{ha_changes_2022}
K.~Ha, \href{https://onlinelibrary.wiley.com/doi/10.1111/phn.12988}{Changes in awareness on face mask use in {Korea}}, Public Health Nursing 39~(2) (2022) 506--508.
\newblock \href {https://doi.org/10.1111/phn.12988} {\path{doi:10.1111/phn.12988}}.
\newline\urlprefix\url{https://onlinelibrary.wiley.com/doi/10.1111/phn.12988}

\bibitem{Brauer1993}
F.~Brauer, C.~Castillo-Chavez, J.~Velasco-Hernandez, C.~U.~B. Unit, C.~U. D.~o. Biometrics, C.~U. D. o. B. S. a.~C. Biology, \href{https://ecommons.cornell.edu/handle/1813/31802}{Recruitment {Effects} in {Heterosexually} {Transmitted} {Disease} {Models}}, eCommons (Aug. 1993).
\newline\urlprefix\url{https://ecommons.cornell.edu/handle/1813/31802}

\bibitem{Brauer2005}
F.~Brauer, \href{https://doi.org/10.1016/j.mbs.2005.07.006}{The kermack{\textendash}{McKendrick} epidemic model revisited}, Mathematical Biosciences 198~(2) (2005) 119--131.
\newblock \href {https://doi.org/10.1016/j.mbs.2005.07.006} {\path{doi:10.1016/j.mbs.2005.07.006}}.
\newline\urlprefix\url{https://doi.org/10.1016/j.mbs.2005.07.006}

\bibitem{Ajbar2021}
A.~Ajbar, R.~T. Alqahtani, M.~Boumaza, \href{https://doi.org/10.3389/fphy.2021.634251}{Dynamics of an {SIR}-based {COVID}-19 model with linear incidence rate, nonlinear removal rate, and public awareness}, Frontiers in Physics 9 (May 2021).
\newblock \href {https://doi.org/10.3389/fphy.2021.634251} {\path{doi:10.3389/fphy.2021.634251}}.
\newline\urlprefix\url{https://doi.org/10.3389/fphy.2021.634251}

\bibitem{hethcote_gonorrhea_1984}
H.~W. Hethcote, J.~A. Yorke, \href{http://link.springer.com/10.1007/978-3-662-07544-9}{Gonorrhea {Transmission} {Dynamics} and {Control}}, Vol.~56 of Lecture {Notes} in {Biomathematics}, Springer Berlin Heidelberg, Berlin, Heidelberg, 1984.
\newblock \href {https://doi.org/10.1007/978-3-662-07544-9} {\path{doi:10.1007/978-3-662-07544-9}}.
\newline\urlprefix\url{http://link.springer.com/10.1007/978-3-662-07544-9}

\bibitem{fentahun_risky_2014}
N.~Fentahun, A.~Mamo, \href{http://www.ajol.info/index.php/ejhs/article/view/102630}{Risky {Sexual} {Behaviors} and {Associated} {Factors} among {Male} and {Female} {Students} in {Jimma} {Zone} {Preparatory} {Schools}, {South} {West} {Ethiopia}: {Comparative} {Study}}, Ethiopian Journal of Health Sciences 24~(1) (2014) 59.
\newblock \href {https://doi.org/10.4314/ejhs.v24i1.8} {\path{doi:10.4314/ejhs.v24i1.8}}.
\newline\urlprefix\url{http://www.ajol.info/index.php/ejhs/article/view/102630}

\bibitem{lucas_schooling_2019}
A.~M. Lucas, N.~L. Wilson, \href{https://www.tandfonline.com/doi/full/10.1080/00220388.2018.1493195}{Schooling, {Wealth}, {Risky} {Sexual} {Behaviour}, and {HIV}/{AIDS} in {Sub}-{Saharan} {Africa}}, The Journal of Development Studies 55~(10) (2019) 2177--2192.
\newblock \href {https://doi.org/10.1080/00220388.2018.1493195} {\path{doi:10.1080/00220388.2018.1493195}}.
\newline\urlprefix\url{https://www.tandfonline.com/doi/full/10.1080/00220388.2018.1493195}

\bibitem{bhattacharyya_modelling_2021}
R.~Bhattacharyya, P.~Konar, \href{https://link.springer.com/10.1007/s40435-020-00692-1}{Modelling the influence of progressive social awareness, lockdown and anthropogenic migration on the dynamics of an epidemic}, International Journal of Dynamics and Control 9~(2) (2021) 797--806.
\newblock \href {https://doi.org/10.1007/s40435-020-00692-1} {\path{doi:10.1007/s40435-020-00692-1}}.
\newline\urlprefix\url{https://link.springer.com/10.1007/s40435-020-00692-1}

\bibitem{Bai2021}
F.~Bai, F.~Brauer, \href{https://www.mdpi.com/2673-3986/2/1/7}{The effect of face mask use on covid-19 models}, Epidemiologia 2~(1) (2021) 75--83.
\newblock \href {https://doi.org/10.3390/epidemiologia2010007} {\path{doi:10.3390/epidemiologia2010007}}.
\newline\urlprefix\url{https://www.mdpi.com/2673-3986/2/1/7}

\bibitem{ghosh_mathematical_2022}
J.~K. Ghosh, S.~K. Biswas, S.~Sarkar, U.~Ghosh, \href{https://linkinghub.elsevier.com/retrieve/pii/S0378475421004043}{Mathematical modelling of {COVID}-19: {A} case study of {Italy}}, Mathematics and Computers in Simulation 194 (2022) 1--18.
\newblock \href {https://doi.org/10.1016/j.matcom.2021.11.008} {\path{doi:10.1016/j.matcom.2021.11.008}}.
\newline\urlprefix\url{https://linkinghub.elsevier.com/retrieve/pii/S0378475421004043}

\bibitem{Brauer2011}
F.~Brauer, \href{https://doi.org/10.1186/1471-2458-11-S1-S3}{A simple model for behavior change in epidemics}, BMC Public Health 11~(1) (2011) S3.
\newblock \href {https://doi.org/10.1186/1471-2458-11-S1-S3} {\path{doi:10.1186/1471-2458-11-S1-S3}}.
\newline\urlprefix\url{https://doi.org/10.1186/1471-2458-11-S1-S3}

\bibitem{Wang2004}
W.~Wang, S.~Ruan, \href{https://www.sciencedirect.com/science/article/pii/S0022247X03008837}{Bifurcations in an epidemic model with constant removal rate of the infectives}, Journal of Mathematical Analysis and Applications 291~(2) (2004) 775--793.
\newblock \href {https://doi.org/https://doi.org/10.1016/j.jmaa.2003.11.043} {\path{doi:https://doi.org/10.1016/j.jmaa.2003.11.043}}.
\newline\urlprefix\url{https://www.sciencedirect.com/science/article/pii/S0022247X03008837}

\bibitem{Alexander2005}
M.~E. Alexander, S.~M. Moghadas, \href{http://www.jstor.org/stable/4096153}{Bifurcation analysis of an sirs epidemic model with generalized incidence}, SIAM Journal on Applied Mathematics 65~(5) (2005) 1794--1816.
\newline\urlprefix\url{http://www.jstor.org/stable/4096153}

\bibitem{Moghadas2006}
S.~M. Moghadas, M.~E. Alexander, \href{https://doi.org/10.1093/imammb/dql011}{{Bifurcations of an epidemic model with non-linear incidence and infection-dependent removal rate}}, Mathematical Medicine and Biology: A Journal of the IMA 23~(3) (2006) 231--254.
\newblock \href {http://arxiv.org/abs/https://academic.oup.com/imammb/article-pdf/23/3/231/1956090/dql011.pdf} {\path{arXiv:https://academic.oup.com/imammb/article-pdf/23/3/231/1956090/dql011.pdf}}, \href {https://doi.org/10.1093/imammb/dql011} {\path{doi:10.1093/imammb/dql011}}.
\newline\urlprefix\url{https://doi.org/10.1093/imammb/dql011}

\bibitem{Xiao2015}
Y.~Xiao, W.~Zhang, G.~Deng, Z.~Liu, \href{https://doi.org/10.1155/2015/745732}{Stability and bogdanov-takens bifurcation of an sis epidemic model with saturated treatment function}, Mathematical Problems in Engineering 2015 (2015).
\newblock \href {https://doi.org/10.1155/2015/745732} {\path{doi:10.1155/2015/745732}}.
\newline\urlprefix\url{https://doi.org/10.1155/2015/745732}

\bibitem{Li2015}
C.~Li, J.~Li, Z.~Ma, Codimension 3 bt bifurcations in an epidemic model with a nonlinear incidence, Discrete Contin. Dyn. Syst. Ser. B 20~(4) (2015) 1107--1116.

\bibitem{Lu2020}
M.~Lu, C.~Xiang, J.~Huang, Bogdanov-takens bifurcation in a sirs epidemic model with a generalized nonmonotone incidence rate, Discrete and Continuous Dynamical Systems-S 13~(11) (2020) 3125--3138.

\bibitem{Zhang2022}
F.~Zhang, W.~Cui, Y.~Dai, Y.~Zhao, Bifurcations of an sirs epidemic model with a general saturated incidence rate, Mathematical Biosciences and Engineering 19~(11) (2022) 10710--10730.

\bibitem{Shan2014}
C.~Shan, H.~Zhu, \href{https://www.sciencedirect.com/science/article/pii/S0022039614002228}{Bifurcations and complex dynamics of an sir model with the impact of the number of hospital beds}, Journal of Differential Equations 257~(5) (2014) 1662--1688.
\newblock \href {https://doi.org/https://doi.org/10.1016/j.jde.2014.05.030} {\path{doi:https://doi.org/10.1016/j.jde.2014.05.030}}.
\newline\urlprefix\url{https://www.sciencedirect.com/science/article/pii/S0022039614002228}

\bibitem{Misra2022}
A.~K. Misra, J.~Maurya, M.~Sajid, \href{https://www.aimspress.com/article/doi/10.3934/mbe.2022541}{Modeling the effect of time delay in the increment of number of hospital beds to control an infectious disease}, Mathematical Biosciences and Engineering 19~(11) (2022) 11628--11656.
\newblock \href {https://doi.org/10.3934/mbe.2022541} {\path{doi:10.3934/mbe.2022541}}.
\newline\urlprefix\url{https://www.aimspress.com/article/doi/10.3934/mbe.2022541}

\bibitem{Liu2016}
Z.~Liu, P.~Magal, D.~Xiao, Bogdanov--takens bifurcation in a predator--prey model, Zeitschrift f{\"u}r angewandte Mathematik und Physik 67~(6) (2016) 137.

\bibitem{van2002reproduction}
P.~Van~den Driessche, J.~Watmough, Reproduction numbers and sub-threshold endemic equilibria for compartmental models of disease transmission, Mathematical biosciences 180~(1-2) (2002) 29--48.

\bibitem{castro2020bifurcation}
J.~A. Castro, F.~Verduzco, Bifurcation analysis in planar quadratic differential systems with boundary, International Journal of Bifurcation and Chaos 30~(07) (2020) 2030017.

\bibitem{kuznetsov2005practical}
Y.~A. Kuznetsov, Practical computation of normal forms on center manifolds at degenerate bogdanov--takens bifurcations, International Journal of Bifurcation and Chaos 15~(11) (2005) 3535--3546.

\end{thebibliography}

\end{document}